\documentclass[11pt]{article}
\usepackage{a4wide}

\usepackage[numbers]{natbib}  

\usepackage[fleqn]{amsmath}
\usepackage{amssymb}
\usepackage{tabularx}
\usepackage{ltablex}
\usepackage{csquotes}
\usepackage{graphicx}
\usepackage{placeins}
\usepackage{eurosym}
\keepXColumns
\usepackage[dvipsnames]{xcolor}
\usepackage{array} 

\newcolumntype{C}[1]{>{\centering\arraybackslash}m{#1}}

\title{Vintage-Based Formulations in Multi-Year Investment Modelling for Energy Systems}
\author{
Ni Wang\textsuperscript{1}, Germán Morales-España \textsuperscript{1,2} \\
\parbox{\textwidth}{
\centering
\textit{\small
\textsuperscript{1}Energy \& Materials Transition, TNO, The Netherlands \\
\textsuperscript{2}Faculty of Electrical Engineering, Mathematics and Computer Science, TU Delft, The Netherlands}
}
}
\date{}

\begin{document}

\maketitle

\begin{abstract}
This paper reviews two established formulations for modelling multi-year energy investments: the simple method, which aggregates all capacity regardless of commissioning year, and the vintage method, which explicitly tracks investments by year to capture differences in technical parameters over time. While the vintage method improves modelling fidelity, it significantly increases model size. To address this, we propose a novel compact formulation that maintains the ability to represent year-specific characteristics while reducing the dimensionality of the model. The proposed compact formulation is implemented in the open-source model \textit{TulipaEnergyModel.jl} and offers a tractable alternative for detailed long-term energy system planning.
\end{abstract}

\section{Introduction}
This document reviews existing formulations for multi-year investment planning in the literature and introduces a novel formulation implemented in the \textit{TulipaEnergyModel.jl} \cite{Tejada-Arango2023a}, focusing on the structural modeling of investments over time, with an emphasis on tracking vintages.
We present the formulations in order of increasing complexity, culminating in a compact representation that advances the state-of-the-art. The primary objective is to reduce model size through acceptable simplifications, particularly regarding capacity constraints and conversion efficiencies.
To maintain focus on the core modeling aspects, we omit the economic representation, which are described in a separate document \cite{Wang2025}. The terminology used here aligns with that of the \textit{TulipaEnergyModel.jl}; readers are encouraged to consult its documentation \cite{Tejada-Arango2023b} for clarification where needed.

\section{Mathematical formulations}

\subsection{Simple investment method}

We now present a simplified formulation of the objective function and the capacity constraints.

For clarity, the objective function omits discounting and several techno-economic cost components that are typically included in more detailed formulations. In addition to these standard terms, other relevant cost categories may also be considered to improve model realism and policy relevance. These include indirect costs \cite{Wang2023}, as well as market-related costs that arise from evolving market designs and trading mechanisms \cite{Estanqueiro2022}. While these are not included in the current formulation, they can be incorporated when needed, particularly in studies aiming for high-fidelity system representation.

\begin{align}
\min & \sum_{a,y} p_{a,y}^{\text{investment cost}} p_{a}^{\mathrm{unit\,capacity}} v^{\text{inv}}_{a, y} \nonumber \\  
& + \sum_{a,y} p_{a,y}^{\text{fixed cost}} p_{a}^{\mathrm{unit\,capacity}} v^{\text{accumulated simple}}_{a, y} \nonumber \\
&+ \sum_{f,y, k_y, b_{k_y}} p^{\text{rp weight}}_{k_y} \cdot p^{\text{variable cost}}_{f, y} \cdot p^{\text{duration}}_{b_{k_y}} \cdot v^{\text{flow}}_{f, k_y, b_{k_y}}
\end{align}

\begin{align}
v^{\text{available units simple method}}_{a,y} & = p^{\text{initial units}}_{a,y} + \sum_{i \in \{\mathcal{Y}^\text{i}_a: y - p^{\text{technical lifetime}}_{a} + 1  \le i \le y \}}  v^{\text{inv}}_{a,i} \nonumber \\
& - \sum_{i \in \{\mathcal{Y}_: y - p^{\text{technical lifetime}}_{a} + 1  \le i \le y \}} v^{\text{decom simple}}_{a,i} \nonumber \\
& \forall a \in \mathcal{A}^{\text{simple investment}} \cup \mathcal{A}^{\text{operation}}, \forall y \in \mathcal{Y} \\
\sum_{f \in \mathcal{F}^{\text{out}}_{a, y}} v^{\text{flow}}_{f,k_y,b_{k_y}} & \leq p^{\text{availability profile}}_{a,k_y,b_{k_y}} \cdot p^{\text{capacity}}_{a, y} \cdot v^{\text{available units simple method}}_{a,y}   \quad  \nonumber 
\\ & \forall y \in \mathcal{Y}, \forall a \in (\mathcal{A}^{\text{simple investment}} \cup \mathcal{A}^{\text{operation}}) \cap \mathcal{A}_y^{\text{p}}, \nonumber \\
& \forall k_y \in \mathcal{K}_y,\forall b_{k_y} \in \mathcal{B}_{k_y} 
\end{align}

This constraint for $v^{\text{available units simple method}}_{a,y}$ tracks the total available units of an asset $a$ in year $y$ using a simplified accumulation method. It is composed of:
\begin{itemize}
    \item The initial capacity available at the beginning of the planning period.
    \item Investments made in previous years $i$ that are still operational in year $y$, i.e., within the technical lifetime of the asset.
    \item Decommissioning possibilities which remove capacity from the available pool.
\end{itemize}

Additional notes:
\begin{itemize}
    \item The capacity constraint applies to producer assets, i.e., technologies that generate usable outputs. Similar constraints can be formulated for other asset types (e.g., storage, conversion) with appropriate adjustments.
    \item Assets are indexed by name and year, enabling investment tracking over time. For instance, the asset name might be wind, with instances in 2030 and 2040 denoted by different indices. 
    \item  The parameter $p^{\text{initial unit}}_{a,y}$ captures preexisting capacity in year $y$, while $v^{\text{inv}}_{a, i}$ represents the capacity added in year $i$ that remains operational in year $y$. The technical lifetime parameter defines the cutoff beyond which assets are no longer active unless decommissioned earlier via $v^{\text{decom simple}}_{a,i}$.
\end{itemize}

\subsection{Vintage investment method}

While this approach simplifies the model and reduces its computational burden, it comes with a notable limitation: it cannot distinguish between units commissioned in different years. Consequently, it is not possible to account for variations in technical characteristics such as conversion efficiencies or availability profiles that may differ across vintages.

To overcome this limitation, the vintage investment method is often employed as done in the national energy system model OPERA \cite{vanStralen2021}. This approach tracks both the year of investment and the year of operation, allowing the model to represent vintage-specific parameters. Implementing this method requires two key changes:

\begin{itemize}
    \item The production variable must include an additional index for the investment year, so that, for any operational year $y$, the contribution of capacity commissioned in each previous year $v$ can be separately modeled.
    \item The capacity constraints must also be defined at the level of individual vintages. This involves indexing the decommissioning variable by both investment and operational year to determine the active capacity available in each period.
\end{itemize}

To illustrate, consider a case where wind power investments are made in years 2030, 2040, and 2050:
\begin{itemize}
    \item Under the simple method, the model includes:
    \begin{itemize}
        \item 3 investment variables (one per year),
        \item 3 decommissioning variables,
        \item 3 production variables (one per year), and
        \item 3 associated capacity constraints.
    \end{itemize}
    \item In contrast, under the vintage method, the model differentiates between investment and operational years. This leads to:
    \begin{itemize}
        \item 3 investment variables,
        \item 6 decommissioning variables (since each operational year must account for multiple vintage contributions),
        \item 6 production variables (one for each combination of vintage and operational year), and
        \item 6 capacity constraints (again, one per vintage-year pair where the asset is active).
    \end{itemize}
\end{itemize}

The vintage method provides greater modeling fidelity at the cost of increased computational complexity \cite{Poncelet2016}. For brevity, we omit the full mathematical formulation here, though it follows naturally from the principles outlined above.

\subsection{Compact investment method}

Since our primary interest lies in capturing the differences in renewable availability profiles across operational years $y$, rather than explicitly modeling vintage-specific production variables, the capacity constraints can be reformulated to eliminate the vintage index $v$. By aggregating the contribution of all relevant investment years directly within the constraint, we retain the ability to apply year-specific availability profiles while significantly reducing the number of production variables and associated constraints. This results in a more compact formulation that balances model accuracy and computational efficiency.

We now present a compact investment formulation, which offers a middle ground between the simple and vintage approaches. It retains the ability to represent year-specific capacity constraints using availability profiles, while avoiding the full complexity of vintage-indexed production variables.

\begin{align}
\min & \sum_{a,y} p_{a,y}^{\text{investment cost}} p_{a}^{\mathrm{unit\,capacity}} v^{\text{inv}}_{a, y} \nonumber \\  
&+ \sum_{(a,y,v) \in \mathcal{D}^{\text{compact investment}}}  p_{a,y,v}^{\text{fixed cost}}      p_{a}^{\mathrm{unit\,capacity}}v^{\text{accumulated compact}}_{a,y,v}   \nonumber \\
&+ \sum_{f,y, k_y, b_{k_y}} p^{\text{rp weight}}_{k_y} \cdot p^{\text{variable cost}}_{f, y} \cdot p^{\text{duration}}_{b_{k_y}} \cdot v^{\text{flow}}_{f, k_y, b_{k_y}}
\end{align}

\begin{align}
v^{\text{available units compact method}}_{a,y,v} & = p^{\text{initial units}}_{a,y,v} + v^{\text{inv}}_{a,v} \nonumber \\
& - \sum_{i \in \{\mathcal{Y}: v < i \le y\} | (a,i,v) \in \mathcal{D}^{\text{compact investment}}} v^{\text{decom compact}}_{a,i,v} \nonumber \\
& \forall a \in \mathcal{A}^{\text{compact investment}}, \forall y \in \mathcal{Y} \\
\sum_{f \in \mathcal{F}^{\text{out}}_{a, y}} v^{\text{flow}}_{f,k_y,b_{k_y}} & \leq \sum_{v \in \mathcal{V} | (a,y,v) \in \mathcal{D}^{\text{compact investment}}} p^{\text{availability profile}}_{a,k_y,b_{k_y}}  \nonumber \\
&\cdot p^{\text{capacity}}_{a, y} \cdot v^{\text{available units compact method}}_{a, y, v}   \quad  \nonumber \\ 
& \forall y \in \mathcal{Y}, \forall a \in \mathcal{A}^{\text{compact investment}} \cap \mathcal{A}_y^{\text{p}}, \nonumber\\
& \forall k_y \in \mathcal{K}_y,\forall b_{k_y} \in \mathcal{B}_{k_y}
\end{align}

This capacity constraint aggregates the capacity contributions from all vintages $v$ that are still operational in year $y$, applying a single availability profile per asset per year. By removing the vintage index from the production variable, we significantly reduce model size.

\subsubsection{Implications of omitting vintage from production variables}
While this compact formulation allows for heterogeneous availability over operational years, it introduces some limitations by removing vintage-specific detail from production modeling:

\begin{itemize}
    \item Efficiencies of producer assets (e.g., thermal units) can still be approximated through the cost parameter $p_{a,y,v}^{\text{variable}}$, which may vary across years. In practice, however, thermal efficiencies are often assumed to remain constant across the time horizon, so this simplification is usually acceptable.
    \item For conversion assets such as electrolyzers, this formulation limits us to a single efficiency value per operational year, regardless of when the unit was built. In cases where technology improves significantly over time, this simplification may result in a conservative estimate. For instance, applying an average or early-year efficiency to all vintages may underestimate hydrogen production and yield a sub-optimal system cost.
\end{itemize}

This compact formulation is particularly valuable when computational tractability is prioritized over modeling precision. It offers a scalable alternative to vintage modeling, especially in large-scale systems where computational resources are limiting factors.

\section{Conclusion}

The choice between simple, compact, and vintage investment formulations presents a trade-off between model accuracy and computational efficiency. The simple method offers ease of implementation and reduced dimensionality but lacks the ability to capture vintage-specific behavior. The vintage method, while the most detailed, significantly increases model size and complexity. The compact method proposed here strikes a balance by allowing for year-specific availability and capacity accounting without introducing vintage-specific production variables. This approach is particularly well-suited for large-scale energy system models where computational efficiency is a priority, and moderate simplifications in technology representation are acceptable.

\bibliographystyle{unsrtnat}      
\bibliography{references}

\end{document}